\documentclass[10pt]{amsart}
\usepackage{epsfig}
\usepackage{amsmath, amssymb, euscript,  enumerate}
\usepackage{bbm}
\usepackage{color}

\newcommand{\supp}{\text{\rm supp}}

\newcommand{\ap}{\alpha}             
\newcommand{\bt}{\beta}
\newcommand{\gm}{\gamma}             \newcommand{\Gm}{\Gamma}
\newcommand{\dt}{\delta}             \newcommand{\Dt}{\Delta}
\newcommand{\vep}{\varepsilon}

\newcommand{\ld}{\lambda}            \newcommand{\Ld}{\Lambda}
             
\newcommand{\vp}{\varphi}
\newcommand{\om}{\omega}             \newcommand{\Om}{\Omega}
\newcommand{\vr}{\varrho}            \newcommand{\iy}{\infty}
            
\newcommand{\f}{\frac}             \newcommand{\el}{\ell}

\newcommand{\fG}{{\mathfrak G}}

\newcommand{\fL}{{\mathfrak L}}
\newcommand{\fM}{{\mathfrak M}}

\newcommand{\fe}{{\mathfrak e}}

\newcommand{\fm}{{\mathfrak m}}

\newcommand{\BN}{{\mathbb N}}

\newcommand{\BR}{{\mathbb R}}

\newcommand{\cE}{{\mathcal E}}
\newcommand{\cF}{{\mathcal F}}

\newcommand{\cI}{{\mathcal I}}

\newcommand{\cK}{{\mathcal K}}

\newcommand{\cO}{{\mathcal O}}

\newcommand{\cS}{{\mathcal S}}
\newcommand{\cT}{{\mathcal T}}

\newcommand{\la}{\langle}          \newcommand{\ra}{\rangle}
\newcommand{\s}{\setminus}         
            \newcommand{\e}{\eta}
\newcommand{\pa}{\partial}        
       
    \newcommand{\ds}{\displaystyle}

 \newcommand{\pf }{\noindent{\it Proof. }}

\newcommand{\aee }{\text{\rm a.e.}} 
  \newcommand{\pv }{\text{\rm p.v.}}
  \newcommand{\dd }{\text{\rm d}}

\newcommand{\rL }{{\text{\rm L}}}

\newcommand{\rX}{{\text{\rm X}}}  \newcommand{\rY}{{\text{\rm Y}}}

\newcommand{\loc}{{\text{\rm loc}}}

\newcommand{\rA }{{\text{\rm A}}}

\newcommand{\rh }{{\text{\rm RH}}}

\newtheorem{thm}[subsection]{Theorem}
\newtheorem{lemma}[subsection]{Lemma}
\newtheorem{cor}[subsection]{Corollary}

\newtheorem{prop}[subsection]{Proposition}

\newtheorem{defn}[subsection]{Definition}

\numberwithin{equation}{section}

\title[ nonlocal Schr\"odinger operators with certain potential]{$L^p$ mapping properties for nonlocal Schr\"odinger operators with certain potential}
\author{ Woocheol Choi and Yong-Cheol Kim  }
\begin{document}
\begin{abstract} In this paper, we consider nonlocal Schr\"odinger equations with certain potentials $V$ given
by an integro-differential operator $L_K$ as follows;
\begin{equation*}L_K u+V u=f\,\,\text{ in $\BR^n$ }\end{equation*}
where $V\in\rh^q$ for $q>\f{n}{2s}$ and $0<s<1$. We denote the solution of the above equation by $\cS_V f:=u$, which is called {\it the inverse of the nonlocal Schr\"odinger operator $L_K+V$ with potential $V$}; that is, $\cS_V=(L_K+V)^{-1}$. Then we obtain a weak Harnack inequality of weak subsolutions of the nonlocal equation 
\begin{equation}\begin{cases}L_K u+V u=0\,\,&\text{ in $\Om$,}\\
                                \quad u=g\,\,&\text{ in $\BR^n\s\Om$,}
\end{cases}\end{equation}
where $g\in H^s(\BR^n)$ and $\Om$ is a bounded open domain in $\BR^n$ with Lipschitz boundary, and also get an improved decay of a fundamental solution $\fe_V$ for $L_K+V$.
Moreover, we obtain $L^p$ and $L^p-L^q$ mapping properties of the inverse $\cS_V$ of the nonlocal Schr\"odinger operator $L_K+V$.
\end{abstract}
\thanks {2000 Mathematics Subject Classification: 47G20, 45K05,
35J60, 35B65, 35D10 (60J75)
}


\address{$\bullet$ Yong-Cheol Kim : Department of Mathematics Education, Korea University, Seoul 136-701,
Korea $\&$ Department of Mathematics, Korea Institute for Advanced Study, Seoul 130-722, Korea  }

\email{ychkim@korea.ac.kr}

\address{$\bullet$ Woocheol Choi : Department of Mathematics, Korea Institute for Advanced Study, Seoul 130-722, Korea } 

\email{wchoi@kias.re.kr}

\maketitle

\tableofcontents

\section{Introduction}
Let $\Om$ be a bounded open domain in $\BR^n$ with Lipschitz boundary.
Then we introduce integro-differential operators of form 
\begin{equation}L_K u(x)=\f{1}{2}\,\pv\int_{\BR^n}\mu(u,x,y)K(y)\,dy, \,\,x\in\Om,
\end{equation}
where $\mu(u,x,y)=2\,u(x)-u(x+y)-u(x-y)$ and the kernel $K:\BR^n\s\{0\}\to\BR_+$
satisfy the property 
\begin{equation}\f{c_{n,s}\,\ld}{|y|^{n+2s}}\le K(y)=K(-y)\le\f{c_{n,s}\,\Ld}{|y|^{n+2s}},\,\,s\in(0,1),\,0<\ld<\Ld<\iy.
\end{equation}
Set $\fL=\{L_K:K\in\cK\}$ where $\cK$ denotes the family of all kernels $K$ satisfying $(1.2)$. In particular, if $K(y)=c_{n,s}|y|^{-n-2s}$ where
$c_{n,s}$ is the normalization constant comparable to $s(1-s)$
given by
\begin{equation}c_{n,s}\int_{\BR^n}\f{1-\cos(\xi_1)}{|\xi|^{n+2s}}\,d\xi=1,
\end{equation}
then $L_K=(-\Delta)^s$ is the fractional Laplacian and it
is well-known that $$\lim_{s\to 1^-}(-\Delta)^s\,u=-\Delta u$$
for any function $u$ in the Schwartz space $\cS(\BR^n)$.

We focus our attention on the nonlocal Schr\"odinger operator $\rL_V=L_K+V$ with potential $V$; as a matter of fact, we consider the nonlocal Schr\"odinger equation with potential $V$ given by
\begin{equation}\begin{cases}\rL_V u=0\,\,&\text{ in $\Om$,}\\
                                \quad u=g\,\,&\text{ in $\BR^n\s\Om$,}
\end{cases}\end{equation}
where $K\in\cK$ and $V$ is a nonnegative potential with $V\in L^1_{\loc}(\BR^n)$.
Then we are interested in $L^p$-estimates and $L^p-L^q$ estimates for the inverse $\cS_V$ of the nonlocal Schr\"odinger operator with nonnegative potential $V$ to be given in Section 6.

Let $\Om$ be a bounded open domain in $\BR^n$ with {\it Lipschitz boundary}
and let $K\in\cK$. Let $\rX(\Om)$ be the linear function space of all real-valued
Lebesgue measurable functions $v$ on $\BR^n$ such that $v|_\Om\in L^2(\Om)$
and
\begin{equation*}\iint_{\BR^{2n}\s(\Om^c\times\Om^c)}\f{|v(x)-v(y)|^2}{|x-y|^{n+2s}}\,dx\,dy<\iy.
\end{equation*}
Set $\rX_0(\Om)=\{v\in\rX(\Om):v=0\,\,\aee\text{ in $\BR^n\s\Om$ }\}.$
For $g\in H^s(\BR^n)$, we consider the convex set of $H^s(\BR^n)$ given by $\rX_g(\Om)=\{v\in H^s(\BR^n):g-v\in\rX_0(\Om)\}$. 

Let $V\in\rh^q$ for $q>\f{n}{2s}$ and $s\in(0,1)$. Then we say that a function $u\in\rX_g(\Om)$ is a {\it weak solution} of the nonlocal
equation (1.4), 
if it satisfies the weak formulation
$$\iint_{\BR^{2n}\s(\Om^c\times\Om^c)}(u(x)-u(y))(\vp(x)-\vp(y))K(x-y)\,dx\,dy+\la Vu,\vp\ra_{L^2(\Om)}=0$$ for any $\vp\in\rX_0(\Om)$.

\,\,We now state our main theorems as follows.

\begin{thm} Let $V\in\rh^q$ for $q>\f{n}{2s}$ with $s\in(0,1)$ and $n\ge 2$. Then there are constants $\vep, C>0$ depending only on $n,\ld,$ and $s$ such that
\begin{equation*}0\le\fe_V(x-y)\le\f{C}{\Xi\bigl(\vep(1+\f{1}{2}|x-y|\,\fm_V(x))^{\f{s}{d_0+1}}\bigr)\,|x-y|^{n-2s}}
\end{equation*} for any $x,y\in\BR^n$, where $\Xi(x)=\sum_{k=0}^\iy\,x^k/(k!)^{\f{n}{2}+s}$, $d_0>0$ is the constant given in Lemma $3.2$ and $\fm_V$ is the fractional auxiliary function to be given in Section $3$.
\end{thm}

The main step in proving Theorem 1.1 is to obtain an improved version of the weak Harnack inequality for weak solutions of the equation (1.4) as follows. To get this, certain type of {\it Caccioppoli estimates} for weak solutions of the equation (1.4) to be obtained in Section 5 will play an important role.

\begin{thm} Let $s\in(0,1)$ and $x_0\in\Om$.
If $u$ is a nonnegative weak solution of the equation $(1.4)$ in $\Om$, then there are universal constants $\vep,C>0$ depending only on $n,s,\ld,\Ld$ such that
\begin{equation*}\sup_{B_{\f{R}{2}}(x_0)}u\le\f{C}{\Xi\bigl(\vep(1+R\,\fm_V(x_0))^{\f{s}{d_0+1}}\bigr)}
\biggl(\f{1}{R^n}\int_{B_R(x_0)}u^2(y)\,dy\biggr)^{\f{1}{2}}.
\end{equation*} for any $R\in(0,\dd(x_0,\pa\Om))$, where $\Xi(x)=\sum_{k=0}^\iy\,x^k/(k!)^{\f{n}{2}+s}$, $d_0>0$ is the constant in Lemma $3.2$ and $\fm_V$ is the fractional auxiliary function to be given in Section $3$.\end{thm}

In the next, we get mapping properties of $L_K\circ\cS_V$ from $L^p(\BR^n)$ into $L^p(\BR^n)$ for $p\in[1,q]$, whenever $V\in\rh^q$ for $q>\f{n}{2s}$ with $s\in(0,1)$ and $n\ge 2$

\begin{thm} If $V\in\rh^q$ is nonnegative for $q>\f{n}{2s}$ with $s\in(0,1)$ and $n\ge 2$, then there is a universal constant $C=C(n,s,q)>0$ such that
\begin{equation*}\|(L_K\circ\cS_V) f\|_{L^p(\BR^n)}\le C\,\|f\|_{L^p(\BR^n)}
\end{equation*} for any $p$ with $1\le p\le q$, where $\cS_V=(L_K+V)^{-1}$.
\end{thm}

\,\,For any $p,q$ with $1<p\le q<\iy$, $s\in(0,1)$, $\theta\in[0,n)$, and $n\ge 2$, we write 
$$\Gm_{\theta}=\biggl\{\biggl(\f{1}{p},\f{1}{q}\biggr)\in(0,1)\times(0,1):\f{1}{p}-\f{1}{q}=\f{\theta}{n},\,p\le q\biggr\}\,\,\text{ and }\,\,W=\fm_V^{2s-\theta}.$$
Denote by $\Dt_{\theta}=\{(p,q)\in(1,\iy)\times(1,\iy):(1/p,1/q)\in\Gm_{\theta}\}$.  Let us introduce the {\it weak $L^p(\BR^n)$ space} which is denoted by $L^{p,\iy}(\BR^n)$.
For $0<p<\iy$, the space $L^{p,\iy}(\BR^n)$ is the class of all real-valued Lebesgue measurable functions $g$ on $\BR^n$ such that
\begin{equation*}\|g\|_{L^{p,\iy}(\BR^n)}:=\sup_{\gm>0}\gm\,\om_g(\gm)^{1/p}<\iy
\end{equation*} 
where $\om_g(\gm)=|\{y\in\BR^n:|g(y)|>\gm\}|$ for $\gm>0$. In fact, it is a quasi-normed linear space for $0<p<\iy$.
Then we obtain the mapping properties of $M_W\circ\cS_V$ in the following theorem, where $M_W$ is the multiplication operator, i.e. $M_W f=W f$. 

\begin{thm}\label{thm-7} Let $V\in\rh^\tau$ be nonnegative for $\tau>\f{n}{2s}$ with $s\in(0,1)$ and $n\ge 2$.

$(a)$ If $\,(p,q)\in\bigl(\bigcup_{\theta\in[0,2s)}\Dt_{\theta}\bigr)\bigcup\{(\iy,\iy)\}$, then we have that
$$\|(M_W\circ\cS_V) f\|_{L^q(\BR^n)}\le C_a\,\|f\|_{L^p(\BR^n)}$$ with a constant $C_a=C_a(n,s,\lambda,p)>0$.

$(b)$ If $p=1$, then there is a universal constant $C_b=C_b(n,s,\ld)>0$ such that
$$\|(M_W\circ\cS_V) f\|_{L^{q,\iy}(\BR^n)}\le C_b\,\|f\|_{L^1(\BR^n)}$$
for any $q\in[1,\f{n}{n-2s})$.

$(c)$ If $(p,q)\in\Dt_{2s}$, then there is a universal constant $C_c=C_c(n,s,\ld,p)>0$ such that
$$\|(M_W\circ\cS_V) f\|_{L^{q,\iy}(\BR^n)}\le C_c\,\|f\|_{L^p(\BR^n)}$$
for any $q\in(\f{n}{n-2s},\iy)$.
\end{thm}

\begin{figure}[htbp]
\begin{center}
\includegraphics[scale=0.6]{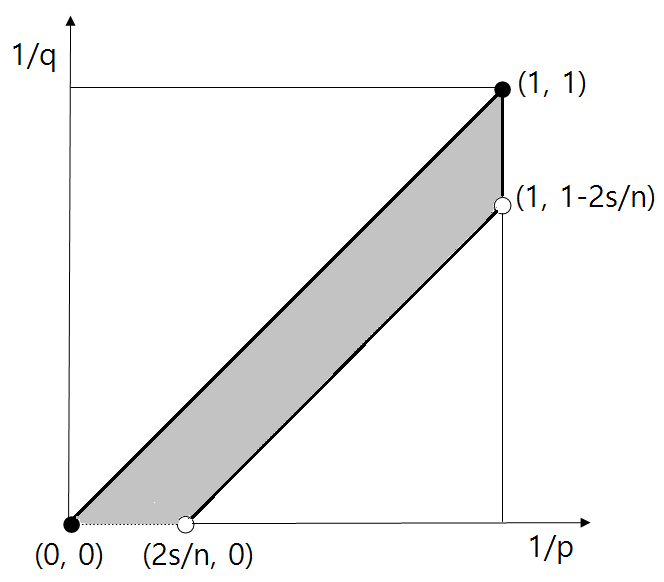}
\caption{The range of $(p,q)$ valid in Theorem \ref{thm-7}} \label{fig:label}
\end{center}
\end{figure}


The paper is organized as follows. In Section 2, we define several function spaces and give the fractional Poincar\'e inequality which was proved in \cite{BBM,MS}.
In Section 3, we introduce the fractional auxiliary function $\fm_V(x)$ related with certain potential $V\in\rh^q$ for $q>\f{n}{2s}$ and $0<s<1$, and deduce a {\it nonlocal version} of the Feffereman-Phong inequality \cite{F} associated with $V$ and $\fm_V$. In Section 4, we also get a weak Harnack inequality for nonnegative weak subsolutions of the equation
\begin{equation}\begin{cases}L_K u=0\,\,&\text{ in $\Om$,}\\
                                \quad u=g\,\,&\text{ in $\BR^n\s\Om$,}
\end{cases}\end{equation}
where $g\in H^s(\BR^n)$ and $\Om$ is a bounded open domain in $\BR^n$ with Lipschitz boundary. In Section 5, 
we furnish a relation between the weak solutions (weak subsolutions and weak supersolutions) of the nonlocal Schr\"odinger equation and the minimizers (subminimizers and superminimizers) of its energy functional, respectively, and also get a {\it Caccioppoli estimate} for weak solutions of the equation. In Section 6, we get an explicit improved upper bound of $\fe_V$. Finally, we obtain $L^p$-estimate and $L^p-L^q$ estimates for the inverse of the nonlocal Schr\"odinger operator $\rL_V$ with potential $V$ in Section 7.

\section{Preliminaries}
Denote by $\cF^n$ the family of all real-valued
Lebesgue measurable functions on $\BR^n$.
Let $\Om$ be a bounded open domain in $\BR^n$ with {\it Lipschitz boundary}
and let $K\in\cK$. Let $\rX(\Om)$ be the linear function space of all
Lebesgue measurable functions $v\in\cF^n$ such that $v|_\Om\in L^2(\Om)$
and
\begin{equation*}\iint_{\BR^{2n}_\Om}\f{|v(x)-v(y)|^2}{|x-y|^{n+2s}}\,dx\,dy<\iy
\end{equation*}
where we denote by $\BR^{2n}_D:=\BR^{2n}\s(D^c\times D^c)$ for a set $D\subset\BR^n$. We
also set
\begin{equation}\rX_0(\Om)=\{v\in\rX(\Om):v=0\,\,\aee\text{ in $\BR^n\s\Om$ }\}
\end{equation}
Since $C^2_0(\Om)\subset\rX_0(\Om)$, we see that $\rX(\Om)$ and $\rX_0(\Om)$ are not
empty. Then we see that $(\rX(\Om),\|\cdot\|_{\rX(\Om)})$ is a normed space,
where the norm $\|\cdot\|_{\rX(\Om)}$ is defined by
\begin{equation}\|v\|_{\rX(\Om)}:=\|v\|_{L^2(\Om)}+\biggl(\iint_{\BR^{2n}_\Om}\f{|v(x)-v(y)|^2}{|x-y|^{n+2s}}\,dx\,dy\biggr)^{1/2}<\iy,
\,\,\,v\in\rX(\Om).
\end{equation}
For $p\ge 1$, let $W^{s,p}(\Om)$ be the usual fractional Sobolev spaces
with the norm
\begin{equation}\|v\|_{W^{s,p}(\Om)}:=\|v\|_{L^p(\Om)}+[v]_{W^{s,p}(\Om)}<\iy
\end{equation}
where the seminorm $[\,\cdot\,]_{W^{s,p}(\Om)}$ is defined by
$$[v]_{W^{s,p}(\Om)}=\biggl(\iint_{\Om\times\Om}
\f{|v(x)-v(y)|^p}{|x-y|^{n+sp}}\,dx\,dy\biggr)^{1/p}.$$
Denote by $H^s(\Om)=W^{s,2}(\Om)$ and let $H^s_c(\BR^n)$ be the class of all functions in $H^s(\BR^n)$ with compact support in $\BR^n$.

By \cite{SV}, there exists a constant $c>1$ depending only on $n,\ld,s$ and $\Om$ such that
\begin{equation}\|v\|_{\rX_0(\Om)}\le\|v\|_{\rX(\Om)}\le c\,\|v\|_{\rX_0(\Om)}
\end{equation}
for any $v\in\rX_0(\Om)$, where
\begin{equation}\|v\|_{\rX_0(\Om)}:=\biggl(\iint_{\BR^{2n}_\Om}\f{|v(x)-v(y)|^2}{|x-y|^{n+2s}}\,dx\,dy\biggr)^{1/2}.
\end{equation} Thus $\|\cdot\|_{\rX_0(\Om)}$ is a norm on $\rX_0(\Om)$ equivalent to (2.4). Moreover it is known \cite{SV} that $(\rX_0(\Om),\|\cdot\|_{\rX_0(\Om)})$ is a
Hilbert space with inner product
\begin{equation}\la
u,v\ra_{\rX_0(\Om)}:=\iint_{\BR^{2n}_\Om}\f{(u(x)-u(y))(v(x)-v(y))}{|x-y|^{n+2s}}\,dx\,dy.
\end{equation}

\begin{lemma} Let $s\in(0,1)$ and $h>0$. If $K\in\cK_0$ and $u\in\rX_0(\Om)$, then we have the following properties; for any $x\in\BR^n$ and $\vr\in(0,h)$,
\begin{equation*}\begin{split}
&(i)\,\,\,\,\vr^{-2}\int_{|x-y|<\vr}|x-y|^2 K(x-y)\,dy+\int_{|x-y|\ge\vr}
K(x-y)\,dy\le\Theta_{n,s}\,\vr^{-2s},\quad\\
&(ii)\,\,\,\,\f{1}{\Ld\,c_{n,s}}\iint_{\Om\times\Om}|u(x)-u(y)|^2 K(x-y)\,dx\,dy
\le\|u\|^2_{W^{s,2}(\Om)}\\
&\qquad\qquad\qquad\qquad\qquad\le\f{1}{\ld\,c_{n,s}}\iint_{\Om\times\Om}|u(x)-u(y)|^2 K(x-y)\,dx\,dy
\end{split}\end{equation*} where $\Theta_{n,s}=\f{\om_n\Ld}{s}$ and $\om_n$ denotes the surface measure of the unit sphere $S^{n-1}$.
\end{lemma}

\pf Refer to \cite{FK} for (i). Also the proof of (ii) is very straightforward. \qed

\,\,\,Next we give the fractional Poincar\'e inequality, which was proved in \cite{BBM,MS}.

\begin{prop} Let $n\ge 1$, $p\ge 1$, $s\in(0,1)$ and $sp<n$. Then there is a universal constant $c_{n,p}>0$ depending only on $n,p$ such that
$$\|u-u_B\|^p_{L^p(B)}\le\f{c_{n,p}(1-s)|B|^{\f{sp}{n}}}{(n-sp)^{p-1}}\,\|u\|^p_{W^{s,p}(B)}$$
for any ball $B\subset\BR^n$.
\end{prop}

\section{The fractional auxiliary function $\fm_V(x)$ }

A locally integrable function in $\BR^n$ that takes values in $[0,\iy)$ almost everywhere is called a {\it weight}. 
For $p\ge 1$ and a weight $w\in L^1_{\rm loc}(\BR^n)$, let $L^p_w(\Om)$ be the weighted $L^p$ class of all real-valued measurable functions $g$ on $\BR^n$ satisfying  
$$\|g\|_{L^p_w(\Om)}:=\biggl(\,\int_{\Om}|g(y)|^p\,w(y)\,dy\biggr)^{\f{1}{p}}<\iy.$$
Consider a class of weights, so-called the {\it Muckenhoupt $\rA_p$-class}, satisfying the following conditions \cite{St}; Let $1\le p<\iy$. Then we say that a weight $w\in L^1_{\loc}(\BR^n)$ satisfies the $\rA_p$-condition (and we denote by $w\in\rA_p$), if
there is a universal constant $C>0$ such that
\begin{equation}\begin{split}
&\biggl(\f{1}{|B|}\int_B w(y)\,dy\biggr)\biggl(\f{1}{|B|}\int_B w(y)^{-\f{1}{p-1}}\biggr)^{p-1}\le C\,\,\,\text{ for $1\le p<\iy$}\\
\end{split}\end{equation} for all balls $B$ in $\BR^n$ (here, the second $L^{\f{1}{p-1}}$-average of $w^{-1}$ on $B$ must  be replaced by $\|w^{-1}\|_{L^{\iy}(B)}$, when $p=1$). The smallest constant $C$ in (3.1) is called the $\rA_p$-norm of $w$ and denoted by $[w]_{\rA_p}$. 
If $w\in\rA_{\iy}:=\cup_{p\ge 1}\rA_p$, then $w\in\rA_p$ for some $p\ge 1$. In this case, it is well-known that $w$ satisfies a reverse H\"older's inequality with exponent $q=1+\e>1$; that is, there are universal constants $C>0$ and $\e>0$ depending only on $n,p$ and $[w]_{\rA_p}$ such that
\begin{equation}\biggl(\f{1}{|B|}\int_B w(y)^q\,dy\biggr)^{\f{1}{q}}
\le\f{C}{|B|}\int_B w(y)\,dy
\end{equation}
for all balls $B$ in $\BR^n$ Let $\rh^q$ be the class of all weights $w$ satisfying (3.2) for some $q>1$, and let $\rh^{\iy}$ be the class of all weights $w$ satisfying that there is a universal constant $C>0$ such that 
$$\|w\|_{L^{\iy}(B)}\le\f{C}{|B|}\int_B w(y)\,dy$$
for all balls $B$ in $\BR^n$. Thus it is obvious that if $w\in\rA_{\iy}$ then there is some $q>1$ such that $w\in\rh^q$, and also $\rh^{\iy}$ is a subclass of any $\rh^q$. Moreover, if $w\in\rh^q$ for $q>1$, then it is well-known that $w\in\rA_{\iy}$, which is equivalent to the following condition; there are some $\ap_0,\bt_0\in(0,1)$ such that
\begin{equation}\bigl|\{x\in B:w(x)\ge\ap_0\,w_B\}\bigr|\ge\bt_0\,|B| 
\end{equation} for all balls $B$ in $\BR^n$, where $w_B=\f{1}{|B|}\int_B w(y)\,dy$ denotes the average of $w$ on $B$.

Throughout this paper, {\it we shall assume that $V\in\rh^q$ for some $q>\f{n}{2s}$ with $s\in(0,1)$ and $n\ge 2$.}
We will consider the auxiliary function to be used in measuring an efficient growth of such weight function $V$, which was introduced by Shen \cite{S}; as a matter of fact, we are considering the nonlocal adaptation of such auxiliary function. 
For $r\in(0,\iy)$, $s\in(0,1)$ and $x\in\BR^n$, we set 
\begin{equation}\fG^s_V(x,r)=\f{1}{r^{n-2s}}\int_{B_r(x)}V(z)\,dz.
\end{equation}
In this case, using H\"older's inequality, it is quite easy to check that
\begin{equation}
\f{\fG^s_V(x,r)}{\fG^s_V(x,R)}\le c_0\biggl(\f{R}{r}\biggr)^{\f{n}{q}-2s}
\end{equation} with a universal constant $c_0>0$, for any $x\in\BR^n$, $s\in(0,1)$ and $0<r<R<\iy$. The assumption $q>\f{n}{2s}$ and (3.5) imply that
\begin{equation}\lim_{r\to 0}\fG^s_V(x,r)=0\,\,\text{ and }\,
\lim_{r\to\iy}\fG^s_V(x,r)=\iy
\end{equation} 
for any $x\in\BR^n$ and $s\in(0,1)$. For $s\in(0,1)$, $x\in\BR^n$ and $V\in\rh^q$ with $q>\f{n}{2s}$, we define
\begin{equation}\f{1}{\fm_V(x)}=\sup\{\rho>0:\fG^s_V(x,\rho)\le 1\}.\end{equation}
From (3.6), it is trivial that $0<\fm_V(x)<\iy$ for all $x\in\BR^n$. If $\fm_V(x)=1/\rho$, then we see that
\begin{equation}\fG^s_V(x,\rho)=1.
\end{equation}
Moreover, it follows from (3.7) that 
\begin{equation}\fG^s_V(x,\rho)\sim 1\,\,\text{ if and only if }\,\,\fm_V(x)\sim\f{1}{\rho}.
\end{equation}
In addition, we mention an well-known fact that the measure $V(z)\,dz$ satisfies the following doubling condition; that is, there is a universal constant $c_1>0$ such that 
\begin{equation}\int_{B_{2r}(x)}V(z)\,dz\le c_1\int_{B_r(x)}V(z)\,dz,
\end{equation}
provided that $V\in\rh^q$ for $q>1$.

\,\,\,For the fractional auxiliary function $\fm_V(x)$, we have the following inequalities whose proof is a nonlocal adaptation of that in \cite{S}.

\begin{lemma} For $s\in(0,1)$, there are universal constants $C_0,d_0>0$ such that 
\begin{equation*}\begin{split}
\text{ $(a)\,\,$ }\fm_V(x)&\sim\fm_V(y)\,\,\text{ if $\,\,|x-y|\le\f{C_0}{\fm_V(x)}$,}\qquad\qquad\qquad\qquad\qquad\qquad\qquad\quad\\
\text{ $(b)\,\,$ }\fm_V(y)&\le C_0\bigl(1+|x-y|\,\fm_V(x)\bigr)^{d_0}\,\fm_V(x)\,\,\text{ for any $x,y\in\BR^n$},\\
\text{ $(c)\,\,$ }\fm_V(y)&\ge\f{\fm_V(x)}{C_0\bigl(1+|x-y|\,\fm_V(x)\bigr)^{\f{d_0}{d_0+1}}}\,\,\,\text{ for any $x,y\in\BR^n$}.
\end{split}\end{equation*} 
\end{lemma}

\pf (a) Assume that $|x-y|\le C_0\rho$ where $\fm_V(x)=1/\rho$. 
Since the measure $V dz$ has the doubling condition (3.10), by (3.8) we see that
$$\fG^s_V(y,\rho)\,\,\sim\,\,
\fG^s_V(x,\rho)=1.$$
Thus it follows from (3.9) that $\fm_V(y)\sim\fm_V(x)$.

(b) Assume that $|x-y|\sim 2^k\rho$ for $k\in\BN$. For $\rho_1\in(0,\rho)$, we choose some $j\in\BN$ so that $2^j\rho_1\sim 2^k\rho$. Then by (3.5) and (3.10) we have that
\begin{equation*}\begin{split}
\int_{B_{\rho_1}(y)}V(z)\,dz&\le c_0\,2^{j[n/q-n]}\int_{B_{2^j\rho_1}(y)}V(z)\,dz
\le C\,2^{j[n/q-n]}\int_{B_{2^k\rho}(y)}V(z)\,dz\\
&\le C\,2^{j[n/q-n]}\int_{B_{2^k\rho}(x)}V(z)\,dz\le c_1^k C\,2^{j[n/q-n]}\int_{B_{\rho}(x)}V(z)\,dz\\
&=c_1^k C\,2^{j[n/q-n]}\rho^{n-2s}.
\end{split}\end{equation*}
Thus this leads us to obtain that
\begin{equation*}\begin{split}\fG^s_V(y,\rho_1)
&\le c_1^k C\,2^{j[n/q-n]}\biggl(\f{\rho}{\rho_1}\biggr)^{n-2s}\le C[2^{n/q-n} c_1]^k\biggl(\f{\rho}{\rho_1}\biggr)^{n/q-2s}.
\end{split}\end{equation*}
If we choose some large number $C_1>0$ so that 
\begin{equation}\begin{split}\fG^s_V(y,\rho_1)
\le a_0\,\,\text{ for $\rho_1\ge C_1^{-k}\rho$, }
\end{split}\end{equation} 
then by (3.7) we conclude that
\begin{equation}\f{1}{\fm_V(y)}\ge C_1^{-k}\rho.\end{equation} 
Hence it follows from (3.11) and (3.12) that
\begin{equation*}\fm_V(y)\le C_1^k\,\fm_V(x)\le C_0\bigl(1+|x-y|\,\fm_V(x)\bigr)^{d_0}\,\fm_V(x)
\end{equation*}
where $d_0=\log_2(C_1)$.

(c) Without loss of generality, we may assume that $|x-y|\ge 1/\fm_V(y)$; for, otherwise it can be shown by (a). From (b), we have that
$$\fm_V(x)\le C_0\bigl(1+|x-y|\,\fm_V(y)\bigr)^{d_0}\,\fm_V(y)\le C_0|x-y|^{d_0}[\fm_V(y)]^{d_0+1}.$$
Therefore we obtain that
$$\fm_V(y)\ge\f{[\fm_V(x)]^{\f{1}{d_0+1}}}{C_0|x-y|^{\f{d_0}{d_0+1}}}\ge\f{\fm_V(x)}{C_0\bigl(1+|x-y|\fm_V(x)\bigr)^{\f{d_0}{d_0+1}}}.$$
Hence we complete the proof. \qed

In the following lemma, we get a {\it nonlocal version} of Feffereman-Phong inequality \cite{F} related with the potential $V$ and the fractional auxiliary function $\fm_V$.

\begin{lemma}  Let $n\ge 1$, $s\in(0,1)$ and $2s<n$. If $u\in H^s_c(\BR^n)$, then there exists a universal constant $C_0=C_0(n,s)>0$ such that
\begin{equation*}\begin{split}
&\int_{\BR^n}|u(x)|^2\,[\fm_V(x)]^{2s}\,dx\le C_0\bigl(\,\|u\|^2_{H^s(\BR^n)}+\|u\|^2_{L^2_V(\BR^n)}\bigr).
\end{split}\end{equation*}
\end{lemma}

\pf We take any $z\in\BR^n$. Let $B=B_{\e}(z)$ and set $\fm_V(z)=1/\e$ for $\e>0$. 
By (3.8), we then observe that
\begin{equation}V_B=\f{1}{|B_1|\e^{2s}}.
\end{equation}
By Proposition 2.2, we have that
\begin{equation}\begin{split}\iint_{B\times B}\f{|u(x)-u(y)|^2}{|x-y|^{n+2s}}\,dx\,dy\ge\f{d_{n,s}}{|B_1|\e^{n+2s}}\iint_{B\times B}|u(x)-u(y)|^2\,dx\,dy
\end{split}\end{equation}
where $d_{n,s}=\f{n-2s}{2\,c_{n,2}(1-s)}$. Also we easily obtain the following equality
\begin{equation}\int_B|u(x)|^2\,V(x)\,dx=\f{c_n}{\e^n}\iint_{B\times B}|u(y)|^2\,V(y)\,dx\,dy
\end{equation} where $c_n=1/|B_1|$.
Since $V\in\rA_{\iy}$, by (3.3) there are universal constants $\ap_0,\bt_0\in(0,1)$ not depending on $B$ such that
\begin{equation}\bigl|\{x\in B:V(x)>\ap_0\,V_B\}\bigr|\ge\bt_0\,|B|.
\end{equation}
By (3.13) and (3.16), adding up (3.14) and (3.15) yields that
\begin{equation}\begin{split}&\iint_{B\times B}\f{|u(x)-u(y)|^2}{|x-y|^{n+2s}}\,dx\,dy
+\int_B|u(y)|^2\,V(y)\,dy\\
&\qquad\ge\f{d_{n,s}\wedge c_n}{2\,\e^n}
\iint_{B\times B}\,\biggl(\f{\ap_0}{|B_1|\e^{2s}}\wedge V(y)\biggr)\,|u(x)|^2\,dx\,dy\\
&\qquad\ge\f{\ap_0\bt_0(d_{n,s}\wedge c_n)}{2\,\e^{2s}}
\int_B|u(x)|^2\,dx.
\end{split}\end{equation}
Thus by (3.17) and (a) of Lemma 3.1 we have that
\begin{equation}\begin{split}
&\int_B|u(x)|^2\,[\fm_V(x)]^{2s}\,dx\le\f{c}{\e^{2s}}\int_B|u(x)|^2\,dx\\
&\qquad\qquad\le C\biggl(\iint_{B\times B}\f{|u(x)-u(y)|^2}{|x-y|^{n+2s}}\,dx\,dy+\int_B|u(x)|^2\,V(x)\,dx\biggr),
\end{split}\end{equation}
where the constant $C=C_{n,s}$ is given by
$$C:=\f{2\,c}{\ap_0\bt_0(d_{n,s}\wedge c_n)}.$$
Applying (a) of Lemma 3.1 again, by (3.18) we obtain that
\begin{equation*}\begin{split}
&\int_B|u(x)|^2\,[\fm_V(x)]^{n+2s}\,dx\\
&\,\le C\biggl(\iint_{B\times B}\f{|u(x)-u(y)|^2}{|x-y|^{n+2s}}[\fm_V(x)]^n\,dx\,dy+\int_B|u(x)|^2\,V(x)\,[\fm_V(x)]^n\,dx\biggr).
\end{split}\end{equation*}
Since $B\times B=\{(x,y)\in\BR^n\times\BR^n:|x-z|\vee|y-z|<\e\}$, 
integrating both sides of the above inequality in $z$ over $\BR^n$ and changing the order of integrations yield that
\begin{equation*}\begin{split}
&\int_{\BR^n}|u(x)|^2\,[\fm_V(x)]^{n+2s}\biggl(\int_{B_{\rho}(x)}dz\biggr)dx\\
&\le C\iint_{\BR^n\times\BR^n}\f{|u(x)-u(y)|^2}{|x-y|^{n+2s}}\,[\fm_V(x)]^n
\biggl(\int_{|z-x|\vee|z-y|<\rho}dz\biggr)dx\,dy\\
&\qquad\qquad\qquad+C\int_{\BR^n}|u(x)|^2\,V(x)\,[\fm_V(x)]^n\biggl(\int_{B_{\rho}(x)}dz\biggr)dx,
\end{split}\end{equation*} where $\rho=1/\fm_V(x)$.
Here we note that
\begin{equation}
\int_{|z-x|\vee|z-y|<\rho}dz\le\int_{B_{\rho}(x)}dz=\f{|B_1|}{[\fm_V(x)]^n}.
\end{equation}
Therefore we complete the proof by using the fact that $u$ is supported in $\Om$. \qed

\begin{lemma} For $V\in\rh^q$ with $q>\f{n}{2s}$ and $s\in(0,1)$, there are some universal constants $d_1>0$ and $C=C(n,s)>0$ such that
$$\f{1}{R^{n-2s}}\int_{B_R(x)}V(y)\,dy\le C\,[R\,\fm_V(x)]^{d_1}\,\,\,\text{ whenever $\,R\,\fm_V(x)\ge 1$. }$$
\end{lemma}

\pf Set $r=1/\fm_V(x)$. If $R\,\fm_V(x)\ge 1$, then we may write $2^{k-1}r\le R<2^k r$ for $k\in\BN$. By (3.8) and the doubling condition (3.10), we have that
\begin{equation*}\int_{B_R(x)}V(y)\,dy\le c_1^k\int_{B_r(x)}V(y)\,dy=c_1^k r^{n-2s}.
\end{equation*}
Thus we conclude that 
\begin{equation*}\f{1}{R^{n-2s}}\int_{B_R(x)}V(y)\,dy\le c_1^k\biggl(\f{r}{R}\biggr)^{n-2s}\le 2^{n-2s}(c_1 2^{2s-n})^k\le C\,[R\,\fm_V(x)]^{d_1},
\end{equation*} where $d_1=\log_2 c_1+2s-n$. Hence we are done. \qed

\section{A weak Harnack inequality }

In this section, we obtain a weak Harnack inequality of the weak solution to the following nonlocal elliptic boundary value problem
\begin{equation}\begin{cases}
L_K u=0 &\text{ in $\Om$,}\\ \,\,\,\,\, u=g &\text{ in $\BR^n\s\Om$,} 
\end{cases}\end{equation} where $g\in H^s(\BR^n)$.
In what follows, we consider a bilinear form by
\begin{equation}\la u,v\ra_K=\iint_{\BR^n\times\BR^n}(u(x)-u(y))(v(x)-v(y))K(x-y)\,dx\,dy\,\,\text{ for $u,v\in\rX(\Om)$.}
\end{equation}
For given $g\in H^s(\BR^n)$, we consider the convex subsets of $H^s(\BR^n)$ by
\begin{equation}\begin{split}&\rX_g^{\pm}(\Om)=\{v\in H^s(\BR^n):(g-v)^{\pm}\in\rX_0(\Om)\},\\
&\rX_g(\Om):=\rX^+_g(\Om)\cap\rX^-_g(\Om)=\{v\in H^s(\BR^n):g-v\in\rX_0(\Om)\}.
\end{split}\end{equation}
The weak formulation of the equation (4.1) is as follows;
if $u\in\rX_g(\Om)$ is a {\it weak solution} of the equation (4.1), then it satisfies that  
\begin{equation}\la
u,\vp\ra_K=0
\end{equation} for all $\vp\in\rX_0(\Om)$. Moreover, we observe that the weak solution $u$ is the minimizer of the energy functional 
\begin{equation}\cE(v)=\iint_{\BR^n\times\BR^n}|v(x)-v(y)|^2\,K(x-y)\,dx\,dy,\,\,v\in\rX_g(\Om).
\end{equation}

\begin{defn} A function $u\in\rX^-_g(\Om)\,(\rX^+_g(\Om))$ is said to be a {\rm weak subsolution (weak supersolution)} of the equation $(4.1)$, if it satisfies that
\begin{equation}\la u,\vp\ra_K\le (\ge)\,0
\end{equation} for every nonnegative $\vp\in\rX_0(\Om)$. Also a function $u$ is a {\rm weak solution} of the equation $(4.1)$, if it is both a weak subsolution and a weak supersolution. So any weak solution $u$ of $(4.1)$ must be in $\rX_g(\Om)$ and satisfy $(4.4)$.
\end{defn}

In the next, we consider the definition of subminimizer and superminimizer of the functional in (4.5) to get better understanding about weak subsolutions and supersolutions.

\begin{defn} Let $g\in H^s(\BR^n)$. $(a)$ A function $u\in\rX^-_g(\Om)$ is said to be a {\rm subminimizer} of the functional $(4.5)$ over $\rX^-_g(\Om)$, if it satisfies that 
$$\,\cE(u)\le\cE(u+\vp)\,\,\,\text{ for all nonpositive $\vp\in\rX_0(\Om)$.}$$ A function $u\in\rX^+_g(\Om)$ is said to be a {\rm superminimizer} of the functional $(4.5)$ over $\rX^+_g(\Om)$, if it satisfies that $$\,\cE(u)\le\cE(u+\vp)\,\,\,\text{ for all nonnegative $\vp\in\rX_0(\Om)$.}$$ 

$(b)$ A function $u$ is said to be a {\rm minimizer} of the functional $(4.5)$ over $\rX_g(\Om)$, if it is both a subminimizer and a superminimizer. So any minimizer $u$ must be in $\rX_g(\Om)$ and satisfies that $\,\cE(u)\le\cE(u+\vp)$ for all $\vp\in\rX_0(\Om)$.
\end{defn}

\begin{thm} If $s\in(0,1)$, then there is a unique minimizer of the functional $(4.5)$. Moreover, a function $u\in\rX^-_g(\Om)\,(\rX^+_g(\Om))$ is a subminimizer $(superminimizer)$ of the functional $(4.5)$ over $\rX^-_g(\Om)\,(\rX^+_g(\Om))$ if and only if it is a weak subsolution $(weak$ $supersolution)$ of the equation $(4.1)$. In particular, a function $u\in\rX_g(\Om)$ is a minimizer of the functional $(4.5)$ if and only if it is a weak solution of the equation $(4.1)$.
\end{thm}

\pf Applying standard way of calculus of variations, we proceed with our proof. 
Take any minimizing sequence $\{u_k\}\subset\rX_g(\Om)$. By the fractional Sobolev inequality \cite{DPV}, we can take a subsequence $\{u_{k_j}\}\subset\rX_g(\Om)$ which converges strongly to $u$ in $L^2(\Om)$. Then there is a subsequence $\{u_{k_i}\}$ of  $\{u_{k_j}\}$ converging $\aee$ in $\Om$ to $u\in\rX_g(\Om)$. Thus it follows from Fatou's lemma that the energy functional $\cE(\cdot)$ is weakly lower semicontinuous on $\rX_g(\Om)$. This implies that $u$ is an actual minimizer of (4.5). Also the uniqueness of the minimizer easily follows from the strict convexity of the functional (4.5). 

Next, we prove the equivalence only for the weak supersolution case, because the other cases can be shown in a similar way. First, if $u\in\rX^+_g(\Om)$, then we note that
\begin{equation}
\cE(u+\vp)-\cE(u)=2\la u,\vp\ra_{\rX(\Om)}+\|\vp\|^2_{\rX_0(\Om)}
\end{equation} for all nonnegative $\vp\in\rX_0(\Om)$.
Thus this implies that a weak supersolution $u\in\rX^+_g(\Om)$ of the equation $(4.1)$ is a superminimizer of the functional $(4.5)$ over $\rX^+_g(\Om)$. On the other hand, we assume that $u\in\rX^+_g(\Om)$ is a superminimizer of the functional $(4.5)$. Then by (4.7) we see that $$2\la u,\vp\ra_{\rX(\Om)}+\|\vp\|^2_{\rX_0(\Om)}\ge 0$$ for all nonnegative $\vp\in\rX_0(\Om)$. 
Since $\vep\vp\in\rX_0(\Om)$ and it is nonnegative for any $\vep>0$ and $\vp\in\rX_0(\Om)$, we have that
$$2\la u,\vp\ra_{\rX(\Om)}+\vep\|\vp\|^2_{\rX_0(\Om)}\ge 0\,\,\,\,\text{ for any $\vep>0$.}$$ Taking $\vep\to 0$, we conclude that $\la u,\vp\ra_{\rX(\Om)}\ge 0$  for any nonnegative $\vp\in\rX_0(\Om)$. Thus $u$ is a weak supersolution of the equation $(4.1)$. Therefore we complete the proof.  
\qed

\,\,As in \cite{DKP2}, we consider the {\it nonlocal tail} $\cT(f;x_0,R)$ of a function $f$ in the open ball $B_R(x_0)\subset\Om$ with center $x_0\in\BR^n$ and radius $R>0$, which plays a crucial role in the regularity results of the nonlocal equations unlike the local equations. For $f\in\rX_0(\Om)$ and $r>0$, we define
\begin{equation*}\cT(f;x_0,r)=r^{2s}\int_{\BR^n\s B_r(x_0)}\f{|f(y)|}{|y-x_0|^{n+2s}}\,dy.
\end{equation*}
We note that this nonlocal tail is well-defined because $\rX_0(\Om)$ is compactly imbedded in $L^2(\Om)$ by the fractional Sobolev inequality on $\rX(\Om)$ \cite{DPV} for a bounded domain $\Om\subset\BR^n$ with Lipschitz boundary.

Next, we mention some results on local boundedness and nonlocal tail properties of weak subsolutions of the equation (4.1) obtained  in \cite{DKP1,DKP2}.

\begin{thm} Let $u\in\rX^-_g(\Om)$ be a weak subsolution of the equation $(4.1)$ where $g\in H^s(\BR^n)$, and let $s\in(0,1)$ and $B_r(x_0)\subset\Om$. Then there is a constant $c_1>0$ depending only on $n,s,\ld$ and $\Ld$ such that
\begin{equation*}\sup_{B_{r/2}(x_0)}u\le\dt\,\cT(u^+;x_0,r/2)+c_1\,\dt^{-\f{n}{4s}}
\biggl(\,\f{1}{|B_r(x_0)|}\int_{B_r(x_0)}[u^+(y)]^2\,dy\biggr)^{1/2}
\end{equation*} for any $\dt\in(0,1]$. Moreover, if $u\ge 0$ in $B_R(x_0)\subset\Om$ with $0<r<R$, then there is a constant $c_2>0$ depending only on $n,s,\ld$ and $\Ld$ such that
\begin{equation*}\cT(u^+;x_0,r)\le c_2\,\sup_{B_r(x_0)}u+c_2\,\biggl(\f{r}{R}\biggr)^{2s}\cT(u^-;x_0,R).
\end{equation*}
\end{thm}

Next we shall prove a weak Harnack inequality of nonnegative weak subsolutions of the equation (4.1) by using Theorem 4.4. Interestingly, this estimate no longer depends on the nonlocal tail term of the weak solution, but the proof is quite simple.

\begin{prop} Let $u\in\rX_g(\Om)$ be a nonnegative weak subsolution of the equation $(1.6)$ where $g\in H^s(\BR^n)$, and let $s\in(0,1)$ and $B_r(x_0)\subset\Om$. Then there is a constant $C>0$ depending only on $n,s,\ld$ and $\Ld$ such that
\begin{equation*}\sup_{B_{r/2}(x_0)}u\le C\,\biggl(\,\f{1}{|B_r(x_0)|}\int_{B_r(x_0)}u^2(y)\,dy\biggr)^{1/2}.
\end{equation*}\end{prop}

\pf We take some $\dt\in(0,1]$ so that $1-\dt c_2>0$ and choose some $R>r$ with $B_R(x_0)\subset\Om$. Then by Theorem 4.4 we have that
\begin{equation*}\begin{split}
\sup_{B_{r/2}(x_0)}u&\le\dt c_2\sup_{B_{r/2}(x_0)}u+\dt c_2\,\biggl(\f{r}{R}\biggr)^{2s}\cT(u^-;x_0,R)\\
&\qquad+c_1\,\dt^{-\f{n}{4s}}
\biggl(\,\f{1}{|B_r(x_0)|}\int_{B_r(x_0)}u^2(y)\,dy\biggr)^{1/2}. 
\end{split}\end{equation*}
Since $\cT(u^-;x_0,R)=0$, we can obtain the required result by taking
$$C=\f{c_1\,\dt^{-\f{n}{4s}}}{1-\dt c_2}.$$
Hence we complete the proof. \qed

\section{Weak solutions and Caccioppoli estimate for nonlocal Schr\"odinger operators $\rL_V$ }
In this section, we give a relation between weak solutions (weak subsolutions and weak supersolutions) for the nonlocal Schr\"odinger equation and minimizers (subminimizers and superminimizers) of its energy functional, respectively. Also, we obtain a certain type of {\it Caccioppoli estimate} for nonnegative weak subsolutions of the nonlocal Schr\"odinger equation.

\begin{defn} Let $g\in H^s(\BR^n)$ and $V\in\rh^q$ for $q>\f{n}{2s}$ and $s\in(0,1)$. Then we say that a function $u\in\rX_g(\Om)$ is a {\rm weak solution} of the nonlocal
equation 
\begin{equation}\begin{cases}\rL_V u=0&\text{ in $\Om$,}\\
u=g&\text{ in $\BR^n\s\Om$, }
\end{cases}\end{equation} 
if it satisfies the weak formulation
\begin{equation}\la u,\vp\ra_K+\la Vu,\vp\ra_{L^2(\BR^n)}=0\,\,\,\text{ for any $\vp\in\rX_0(\Om)$.} \end{equation}
\end{defn}
In fact, it turns out that the weak solution of the equation (5.1) is the minimizer of the energy functional
\begin{equation}\cE_V(v)=\cE(v)+\|v\|^2_{L^2_V(\BR^n)},\,\,v\in\rY_g(\Om):=\rX_g(\Om)\cap L^2_V(\BR^n),
\end{equation} where $g\in H^s(\BR^n)$.
We consider function spaces $\rY^+_g(\Om)$ and $\rY^-_g(\Om)$ defined by $$\rY^{\pm}_g(\Om)=\{v\in\rY_g(\Om):(g-v)^{\pm}\in\rX_0(\Om)\}.$$ Then we see that $\rY_g(\Om)=\rY^+_g(\Om)\cap\rY^-_g(\Om)$. When $u=g=0$ in $\BR^n\s\Om$, we easily see that $\rY_0(\Om)=\rX_0(\Om)\cap L^2_V(\Om)$ and $\rY_0(\Om)$ is a Hilbert space with the inner product defined by
$\la u,v\ra_{\rY_0(\Om)}=\la u,v\ra_{\rX_0(\Om)}+\la V u,v\ra_{L^2(\Om)}.$
Moreover, we see that $\rY_0(\Om)=\rX_0(\Om)$ and they are norm-equivalent (refer to \cite{CK}).
 
As in Section 4, we define weak subsolutions and weak supersolutions of the nonlocal equation (5.1) in the following definition. 

\begin{defn} Let $g\in H^s(\BR^n)$. A function $u\in\rY^-_g(\Om)\,(\rY^+_g(\Om))$ is said to be a {\rm weak subsolution (weak supersolution)} of the equation $(5.1)$, if it satisfies that
\begin{equation}\la u,\vp\ra_K+\la Vu,\vp\ra_{L^2(\BR^n)}\le(\ge)\,\,0
\end{equation} for every nonnegative $\vp\in\rX_0(\Om)$. Also a function $u$ is a {\rm weak solution} of the equation $(5.1)$, if it is both a weak subsolution and a weak supersolution. So any weak solution $u$ of $(5.1)$ must be in $\rY_g(\Om)$ and satisfy $(5.2)$.
\end{defn}

In the next, we furnish the definition of subminimizer and superminimizer of the functional (5.3) to understand well weak subsolutions and supersolutions of the nonlocal Schr\"odinger equation (5.1).

\begin{defn} Let $g\in H^s(\BR^n)$. $(a)$ A function $u\in\rY^-_g(\Om)$ is said to be a {\rm subminimizer} of the functional $(5.3)$ over $\rY^-_g(\Om)$, if it satisfies that
$$\,\cE_V(u)\le\cE_V(u+\vp)\,\,\,\,\text{ for all nonpositive $\vp\in\rX_0(\Om)$.}$$ Also, a function $u\in\rY^+_g(\Om)$ is said to be a {\rm superminimizer} of the functional $(5.3)$ over $\rY^+_g(\Om)$, if it satisfies that $$\,\cE_V(u)\le\cE_V(u+\vp)\,\,\,\,\text{ for all nonnegative $\vp\in\rX_0(\Om)$.}$$

$(b)$ A function $u$ is said to be a {\rm minimizer} of the functional $(5.3)$ over $\rY_g(\Om)$, if it is both a subminimizer and a superminimizer. So any minimizer $u$ must be in $\rY_g(\Om)$ and satisfy that $\,\cE_V(u)\le\cE_V(u+\vp)$ for all $\vp\in\rX_0(\Om)$.
\end{defn}

Let $\rY(\Om)$ be the normed subspace of $\rX(\Om)$ which is endowed with the norm
$$\|u\|_{\rY(\Om)}:=\sqrt{\|u\|^2_{\rX(\Om)}+\|u\|^2_{L^2_V(\Om)}}<\iy,\,\,u\in\rY(\Om).$$
In order to obtain the existence of the minimizer of the functional $\cE_V$ on $\rY_g(\Om)$, we need a compactness theorem $\rY(\Om)\hookrightarrow L^2(\Om)$ as follows.

\begin{thm} Let $n\ge 1$, $s\in(0,1)$ and $2s<n$. If $\,u\in\rY(\Om)$, then there exists a universal constant $C>0$ depending on $n,s$ and $\ld$ such that
$$\|u\|_{L^2(\Om)}\le C\,\|u\|_{\rY(\Om)}.$$
Moreover, any bounded sequence in $\rY(\Om)$ is precompact in $L^2(\Om)$.
\end{thm}

\pf Since $\|u\|_{H^s(\Om)}\le c(\ld)\|u\|_{\rX(\Om)}\le c(\ld)\|u\|_{\rY(\Om)}$ for any $u\in\rY(\Om)$, it follows from the fractional Sobolev inequality \cite{DPV} that
$$\|u\|_{L^2(\Om)}\le C\,\|u\|_{H^s(\Om)}\le C\,\|u\|_{\rY(\Om)}.$$
Thus the precompactness in $L^2(\Om)$ can be obtained by weak compactness. \qed

\begin{lemma} If $s\in(0,1)$, then there is a unique minimizer of the functional $(5.3)$. Moreover, a function $u\in\rY^-_g(\Om)\,(\rY^+_g(\Om))$ is a subminimizer $(superminimizer)$ of the functional $(5.3)$ over $\rY^-_g(\Om)\,(\rY^+_g(\Om))$ if and only if it is a weak subsolution $(weak$ $supersolution)$ of the equation $(5.1)$. In particular, a function $u\in\rY_g(\Om)$ is a minimizer of the functional $(5.3)$ if and only if it is a weak solution of the nonlocal equation $(5.1)$.
\end{lemma}

\pf We proceed with our proof as in Lemma 4.4. Take any minimizing sequence $\{u_k\}\subset\rY_g(\Om)$. By applying Theorem 5.4, we can take a subsequence $\{u_{k_j}\}\subset\rY_g(\Om)$  such that
$$u_{k_j}\to u\,\,\text{ in $L^2(\Om)$}$$ as $j\to\iy$.
So there exist a subsequence $\{u_{k_i}\}$ of $\{u_{k_j}\}$ which converges $\aee$ in $\Om$ to $u\in\rY_g(\Om)$. Thus, by applying Fatou's lemma, we can show that the energy functional $\cE_V(\cdot)$ is weakly semicontinuous in $\rY_g(\Om)$. This implies that $u$ is a minimizer of (5.3). The uniqueness of the minimizer also follows from the strict convexity of the functional (5.3). 

Next, we show the equivalency only for the weak supersolution case, because the other case can be done in a similar way. First, if $u\in\rY^+_g(\Om)$, then we observe that
\begin{equation}\begin{split}&\cE_V(u+\vp)-\cE_V(u)\\
&\qquad\qquad=2\la u,\vp\ra_K+2\la Vu,\vp\ra_{L^2(\BR^n)}+\|\vp\|^2_{\rX_0(\Om)}+\|\vp\|^2_{L^2_V(\Om)}
\end{split}\end{equation} for all nonnegative $\vp\in\rX_0(\Om)$. This implies that a weak supersolution $u\in\rY_g^+(\Om)$ of the equation (5.1) is a superminimizer of the functional (5.3) over $\rY_g^+(\Om)$. 

On the other hand, we suppose that $u\in\rY^+_g(\Om)$ is a superminimizer 
of the functional (5.3). Then it follows from (5.5) that
$$2\la u,\vp\ra_K+2\la Vu,\vp\ra_{L^2(\BR^n)}+\|\vp\|^2_{\rX_0(\Om)}+\|\vp\|^2_{L^2_V(\Om)}\ge 0$$ for all nonnegative $\vp\in\rX_0(\Om)$. Since $\vep\vp\in\rX_0(\Om)$ and it is nonnegative for any $\vep>0$ and $\vp\in\rX_0(\Om)$, we obtain that
$$2\la u,\vp\ra_K+2\la Vu,\vp\ra_{L^2(\BR^n)}+\vep\|\vp\|^2_{\rX_0(\Om)}+\vep\|\vp\|^2_{L^2_V(\Om)}\ge 0$$ for any $\vep>0$. Taking $\vep\to\iy$, we can conclude that
$\la u,\vp\ra_K+\la Vu,\vp\ra_{L^2(\BR^n)}\ge 0$ for any nonnegative $\vp\in\rX_0(\Om)$.
Hence $u$ is a weak supersolution of the equation (5.1). Therefore we are done. \qed

\begin{lemma} If $\ap,\bt\in\BR$ and $a,b\ge 0$, then we have the equality
$$(\bt-\ap)(b^2\bt-a^2\ap)=(b\bt-a\ap)^2-\ap\bt(b-a)^2$$
\end{lemma}

\pf By simple calculation, we have that
\begin{equation*}\begin{split}
(\bt-\ap)(b^2\bt-a^2\ap)&=b^2\bt^2-2ab\ap\bt+a^2\ap^2+2ab\ap\bt-b^2\ap\bt-a^2\ap\bt\\
&=(b\bt-a\ap)^2-\ap\bt(b-a)^2.
\end{split}\end{equation*}
Hence we are done. \qed

\,\,Next we will prove the following type of {\it Caccioppoli estimate} for weak solutions of the equation (5.1).

\begin{lemma} Let $s\in(0,1)$ and $x_0\in\Om$. Suppose that $u$ is a nonnegative weak subsolution of the nonlocal equation $(5.1)$ in $\Om$. Then there is a constant $C>0$ depending only on $n,s,\ld$ and $\Ld$ such that
\begin{equation*}  
\|u\|^2_{L^2_V(B_{R_*}(x_0))}+\ld\,c_{n,s}\,\|\phi\,u\|^2_{H^s(\BR^n)}\le\f{20\,\Theta_{n,s}+C}{(R-r)^{2s}}\biggl(\f{R}{R-r}\biggr)^n\|u\|^2_{L^2(B_R(x_0))}
\end{equation*} for any $r\in(0,\dd(x_0,\pa\Om)/2)$ and any $R\in(r,2r]$, where $\Theta_{n,s}$ is the constant in Lemma 2.3 and $\phi$ is the function defined by
$$\phi(x)=\phi_{r,R_*,x_0}(x):=\biggl(\f{R_*-|x-x_0|}{R_*-r}\vee 0\biggr)
\wedge 1,\text{ where $r<R_*=(R+r)/2<R$.}$$
\end{lemma}

\pf We use $\vp(x)=\phi^2(x)u(x)$ as a testing function in (5.2). Then we note that
\begin{equation}\la u,\vp\ra_{\rX(\Om)}+\la Vu,\vp\ra_{L^2(\Om)}=\la u,\vp\ra_K+\la Vu,\vp\ra_{L^2(\BR^n)}\le 0.
\end{equation}
Applying Lemma 5.6, we obtain that
\begin{equation*}\begin{split}\la u,\vp\ra_{\rX(\Om)}&=\iint_{\BR^{2n}_\Om}
\bigl(u(x)-u(y)\bigr)\bigl(\vp(x)-\vp(y)\bigr)K(x-y)\,dx\,dy\\
&=\iint_{B^2_r(x_0)}\bigl(u(x)-u(y)\bigr)^2\,K(x-y)\,dx\,dy\\
&+\iint_{\BR^{2n}_\Om\s B^2_r(x_0)}\bigl(\phi(x)u(x)-\phi(y)u(y)\bigr)^2\,K(x-y)\,dx\,dy\\
&-\iint_{\BR^{2n}_\Om\s B^2_r(x_0)}\bigl(\phi(x)-\phi(y)\bigr)^2
\,u(x)\,u(y)\,K(x-y)\,dx\,dy
\end{split}\end{equation*}
where $B^2_r(x_0):=B_r(x_0)\times B_r(x_0)$. Thus by (5.6) we have that
\begin{equation}\begin{split}&\int_{B_r(x_0)}V(x)\,u^2(x)\,dx+\iint_{\BR^{2n}_\Om}
\bigl(\phi(x)u(x)-\phi(y)u(y)\bigr)^2\,K(x-y)\,dx\,dy\\
&\quad\le\iint_{\BR^{2n}_\Om\s B^2_r(x_0)}\bigl(\phi(x)-\phi(y)\bigr)^2\,u(x)\,u(y)
\,K(x-y)\,dx\,dy:=\cI.
\end{split}\end{equation}
From the property of $\phi$, we see that 
\begin{equation}\sup_{x,y\in\BR^n}\f{\bigl(\phi(x)-\phi(y)\bigr)^2}{|x-y|^2}
\le\biggl(\f{1}{R_*-r}\biggr)^2\le\f{4}{(R-r)^2}.
\end{equation}
Since we have the estimate
\begin{equation*}|y-x|\ge|y-x_0|-|x-x_0|\ge\f{(R-r)|y-x_0|}{2R}\,\,\text{ for any $(x,y)\in B_{R_*}(x_0)\times B^c_R(x_0)$, }
\end{equation*} 
it follows from Lemma 2.1, (5.7), (5.8) and Cauchy's inequality that
\begin{equation}\begin{split}
\cI&=\f{1}{2}\iint_{B^2_R(x_0)\s B^2_r(x_0)}\bigl(\phi(x)-\phi(y)\bigr)^2\,(u^2(x)+u^2(y))\,K(x-y)\,dx\,dy\\
&\qquad+2\iint_{B_R(x_0)\times B^c_R(x_0)}\phi^2(x)
\,u(x)\,u(y)\,K(x-y)\,dx\,dy\\
&\le\iint_{B^2_R(x_0)}\bigl(\phi(x)-\phi(y)\bigr)^2\,u^2(x)\,K(x-y)\,dx\,dy\\
&\qquad+2\int_{B_R(x_0)}\phi^2(x)
\,u(x)\biggl(\int_{B^c_{2R}(x_0)}\,u(y)\,K(x-y)\,dy\biggr)\,dx\\
&\le\f{4\,\Theta_{n,s}}{(R-r)^{2s}}\,\|u\|^2_{L^2(B_R(x_0))}\\
&\quad+\Ld\,c_{n,s}\biggl(\f{2R}{R-r}\biggr)^{n+2s}\|u\|_{L^1(B_R(x_0))}\int_{B^c_R(x_0)}\f{|u(y)|}{|y-x_0|^{n+2s}}\,dy.
\end{split}\end{equation} 
Since $\cT(u^-;x_0,R)=0$, it follows from Theorem 4.4 and Proposition 4.5 that 
\begin{equation*}\begin{split}
\cI&-\f{4\,\Theta_{n,s}}{(R-r)^{2s}}\,\|u\|^2_{L^2(B_R(x_0))}\\
&\le 2^{n+2s}\Ld\,c_{n,s}\biggl(\f{R}{R-r}\biggr)^{n+2s}|B_R(x_0)|^{\f{1}{2}}\|u\|_{L^2(B_R(x_0))}\biggl(\f{R}{2}\biggr)^{-2s}\cT(u;x_0,R/2)\\
&\le\f{2^{n+4}\,c_2\Ld\,c_{n,s}}{(R-r)^{2s}}\biggl(\f{R}{R-r}\biggr)^n|B_R(x_0)|^{\f{1}{2}}\|u\|_{L^2(B_R(x_0))}\,\sup_{B_{R/2}(x_0)}u\\
&\le\f{C}{(R-r)^{2s}}\biggl(\f{R}{R-r}\biggr)^n|B_R(x_0)|^{\f{1}{2}}\|u\|_{L^2(B_R(x_0))}
\biggl(\,\f{1}{|B_R(x_0)|}\int_{B_R(x_0)}u^2(y)\,dy\biggr)^{\f{1}{2}}\\
&=\f{C}{(R-r)^{2s}}\biggl(\f{R}{R-r}\biggr)^n\|u\|^2_{L^2(B_R(x_0))}.
\end{split}\end{equation*}
Thus, by (5.7), we obtain that
\begin{equation}\begin{split}&\int_{B_r(x_0)}V(x)|u(x)|^2\,dx
+\ld\,c_{n,s}\iint_{\BR^{2n}_\Om}\f{|\phi(x)u(x)-\phi(y)u(y)|^2}{|x-y|^{n+2s}}\,dx\,dy\\
&\qquad\qquad\qquad\qquad\le\cI\le
\f{4\,\Theta_{n,s}+C}{(R-r)^{2s}}\biggl(\f{R}{R-r}\biggr)^n\|u\|^2_{L^2(B_R(x_0))}.
\end{split}\end{equation}
Using $\vp(x)=\phi_0^2(x)u(x)$ as a testing function where $\phi_0=\phi_{R_*,R^*,x_0}$ for $r<R_*<R^*=(3R+r)/4<R,$ via the same way as the above we arrive at the estimate
\begin{equation*}\begin{split}&\int_{B_{R_*}(x_0)}V(x)|u(x)|^2\,dx
+\ld\,c_{n,s}\iint_{\BR^{2n}_\Om}\f{|\phi_0(x)u(x)-\phi_0(y)u(y)|^2}{|x-y|^{n+2s}}\,dx\,dy\\
&\qquad\qquad\qquad\qquad\le\cI\le
\f{16\,\Theta_{n,s}+C}{(R-r)^{2s}}\biggl(\f{R}{R-r}\biggr)^n\|u\|^2_{L^2(B_R(x_0))}.
\end{split}\end{equation*}
Combining this with (5.10), the required estimte can be achieved. \qed

\,\,\,{\bf [Proof of Theorem 1.2]} Let $s\in(0,1)$ and $x_0\in\Om$. Given $k\in\BN$, let $\el_i=\f{3}{4}+\f{i-1}{4k}$ and $R_i=\el_i R$  for $i=1,\cdots,k+1$,
and also let $R^*_i=R_i+\f{1}{8k}R$ for $\,i=1,\cdots,k.\,\,$ For $i=1,\cdots,k$, let us denote by $$\phi_i(x)=\biggl(\f{R^*_i-|x-x_0|}{R^*_i-R_i}\vee 0\biggr)\wedge 1$$ with $R_i<R^*_i=\f{R_i+R_{i+1}}{2}<R_{i+1}$. Fix any $i=1,\cdots,k$ and any $R\in(0,\dd(x_0,\pa\Om))$.
Applying Lemma 3.2 with $\phi_i u$, it follows from Lemma 5.7 that
\begin{equation}\begin{split}\int_{B_{R_i}(x_0)}|u(y)|^2\,[\fm_V(y)]^{2s}\,dy&\le C_{n,s}\bigl(\,\|u\|^2_{L^2_V(B_{R^*_i}(x_0))}+\|\phi_i u\|^2_{H^s(\BR^n)}\,\bigr)\\
&\le\f{C\,k^{n+2s}}{R^{2s}}\,\|u\|^2_{L^2(B_{R_{i+1}}(x_0))},
\end{split}\end{equation} because $R_{i+1}-R_i=\f{1}{4k}R$ and $\f{R_{i+1}}{R_{i+1}-R_i}\le 4k$. 
From (c) of Lemma 3.1 and (5.11), we obtain that
\begin{equation*}
\|u\|^2_{L^2(B_{R_i}(x_0))}\le\f{C\,k^{n+2s}(1+R\fm_V(x_0))^{\f{2s d_0}{d_0+1}}}{R^{2s}\,[\fm_V(x_0)]^{2s}}\,\|u\|^2_{L^2(B_{R_{i+1}}(x_0))}.
\end{equation*}
Continuing this process $k$-times from $i=1$ to $i=k$ yields that
\begin{equation*}\begin{split}
\|u\|^2_{L^2\bigl(B_{\f{3R}{4}}(x_0)\bigr)}\le\f{C^k\,k^{(n+2s)k}(1+R\fm_V(x_0))^{\f{2sk d_0}{d_0+1}}}{R^{2sk}\,[\fm_V(x_0)]^{2sk}}\,\|u\|^2_{L^2(B_R(x_0))}.
\end{split}\end{equation*}
Applying Proposition 4.5, we can derive from the above inequality that
\begin{equation}\sup_{B_{\f{R}{2}}(x_0)}u\le\f{C^{\f{k}{2}}\,k^{(\f{n}{2}+s)k}(1+R\,\fm_V(x_0))^{\f{sk d_0}{d_0+1}}}{R^{sk}\,[\fm_V(x_0)]^{sk}}\biggl(\f{1}{R^n}\int_{B_R(x_0)}u^2(y)\,dy\biggr)^{\f{1}{2}}.
\end{equation}
From the well-known Stirling's formula $k^k\sim e^k k!(2\pi k)^{-1/2}$ as $k\to\iy$, we see that there is a constant $c_0>0$ such that $k^k\le c_0\,e^k k!$ for any $k\in\BN$. 
Combining (5.12) with Proposition 4.5 yields that 
\begin{equation}\begin{split}\f{\sup_{B_{\f{R}{2}}(x_0)}u}{\bigl(\f{1}{R^n}\int_{B_R(x_0)}u^2(y)\,dy\bigr)^{\f{1}{2}}}&\le C^{\f{k}{2}}\,k^{(\f{n}{2}+s)k}\biggl(1\wedge\f{(1+R\,\fm_V(x_0))^{\f{sk d_0}{d_0+1}}}{R^{sk}\,[\fm_V(x_0)]^{sk}}\biggr)\\
&\le\f{c_0^s\,C^{\f{k}{2}}\,e^{(\f{n}{2}+s)k}(k!)^{\f{n}{2}+s}}{(1+R\,\fm_V(x_0))^{\f{sk}{d_0+1}}}\,\,\,\text{ for all $k\in\BN$.}
\end{split}\end{equation} 
Multiplying (5.13) by 
$[(1+R\,\fm_V(x_0))^{\f{sk}{d_0+1}}\vep^k]/(k!)^{\f{n}{2}+s}$ and adding up on $k\in\BN$ where the constant $\vep>0$ is chosen so small that 
$\vep\,e^{\f{n}{2}+s}\,C^{1/2}<1$, we obtain that
\begin{equation*}\begin{split}
&\Xi\bigl(\vep(1+R\,\fm_V(x_0))^{\f{s}{d_0+1}}\bigr)\biggl(\,\sup_{B_{\f{R}{2}}(x_0)}u\biggr)\\
&\qquad\qquad\qquad=\biggl(\,\sup_{B_{\f{R}{2}}(x_0)}u\biggr)\sum_{k=0}^\iy
\f{\bigl(\vep(1+R\,\fm_V(x_0))^{\f{s}{d_0+1}}\bigr)^k}{(k!)^{\f{n}{2}+s}}\\
&\qquad\qquad\qquad\le c^s_0\sum_{k=0}^\iy\bigl(\vep\,e^{\f{n}{2}+s}\,C^{1/2}\bigr)^k\biggl(\f{1}{R^n}\int_{B_R(x_0)}u^2(y)\,dy\biggr)^{\f{1}{2}}\\
&\qquad\qquad\qquad\le C\biggl(\f{1}{R^n}\int_{B_R(x_0)}u^2(y)\,dy\biggr)^{\f{1}{2}}.
\end{split}\end{equation*} Hence we conclude that
$$\quad\qquad\sup_{B_{\f{R}{2}}(x_0)}u\le\f{C}{\Xi\bigl(\vep(1+R\,\fm_V(x_0))^{\f{s}{d_0+1}}\bigr)}
\biggl(\f{1}{R^n}\int_{B_R(x_0)}u^2(y)\,dy\biggr)^{\f{1}{2}}. \quad\qquad\qed$$

\,\,\,{\bf [Proof of Theorem 1.1]} 
Let $\fe_V$ be a fundamental solution for the operator $\rL_V=L_K+V$.
By Theorem 1.1 \cite{CK}, we have that
\begin{equation}0\le\fe_V(x-y)\le\f{C}{|x-y|^{n-2s}}.
\end{equation}
Take any $x\in\BR^n$. Since $\supp(\dt_y)=\{y\}$, we see that $u(z):=\fe_V(z-y)$ satisfies the nonlocal equation (5.1) on $B_R(x)$ where $R=|x-y|/2$. Applying Theorem 1.2 to $u(z)$, we obtain that
\begin{equation*}\begin{split}
\fe_V(x-y)\le\sup_{B_{\f{R}{2}}(x)}|u|\le\f{C}{\Xi\bigl(\vep\bigl(1+R\,\fm_V(x)\bigr)^{\f{s}{d_0+1}}\bigr)}\biggl(\f{1}{R^n}\int_{B_R(x)}|u(z)|^2\,dz\biggr)^{\f{1}{2}}.
\end{split}\end{equation*}
Since $R=|x-y|/2$, we see that $|z-y|\ge|x-y|-|z-x|\ge|x-y|/2$ for any $z\in B_R(x)$, and thus we have that
\begin{equation*}\begin{split}
\fe_V(x-y)\le\f{C}{\Xi\bigl(\vep\bigl(1+\f{1}{2}|x-y|\,\fm_V(x)\bigr)^{\f{s}{d_0+1}}\bigr)\,|x-y|^{n-2s}}.
\end{split}\end{equation*}
Hence we are done. \qed

\begin{cor} Let $V\in\rh^q$ be a nonnegative potential for $q>\f{n}{2s}$ with $s\in(0,1)$ and $n\ge 2$. Then for any $N>0$ there exists a constant $C_N>0$ possibly depending on $n,\ld,s$ such that \begin{equation*}\begin{split}0\le\fe_V(x-y)\le\f{C_N} {\bigl(1+|x-y|\,\fm_V(x)\bigr)^N\,|x-y|^{n-2s}}\,\,\,\,\text{ for $x,y\in\BR^n$.}
\end{split}\end{equation*} 
\end{cor}

\pf For any $N\in(0,\iy)$, it is easy to check that
$$C_N=\sup_{t>0}\,\f{t^N}{\Xi(t)}<\iy.$$ Hence the required estimate immediately follows from Theorem 1.1. \qed

\section{$L^p$ and $L^p-L^q$ mapping properties of the inverse of the nonlocal Schr\"odinger operator }
In this section, we consider the nonhomogeneous nonlocal Schr\"odinger equation with potential $V$ given by
\begin{equation}L_K u+Vu=f\,\,\text{ in $\BR^n$,}
\end{equation} where $V\in\rh^q$ is nonnegative for $q>\f{n}{2s}$ with $s\in(0,1)$ and $n\ge 2$. Then we see that the function
\begin{equation}u(x)=\int_{\BR^n}\fe_V(x-y)f(y)\,dy
\end{equation}
is a solution of the equation (6.1). We denote the solution by $\cS_V f(x):=u(x)$, and so we may write $\cS_V=(L_K+V)^{-1}$. We call $\cS_V$ the {\it inverse of the nonlocal Schr\"odinger operator with nonnegative potentials} $V$.

\,\,\,{\bf [Proof of Theorem 1.3.]} Set $r=1/\fm_V(x)$.  Then by (6.2) we may write
\begin{equation}\begin{split}\cS_V f(x)&=\int_{B_r(x)}\fe_V(x-y)f(y)\,dy+\int_{\BR^n\s B_r(x)}\fe_V(x-y)f(y)\,dy\\
&:=\cS^1_V f(x)+\cS^2_V f(x).
\end{split}\end{equation}
By (5.14) and H\"older's inequality, we have that
\begin{equation}|\cS^1_V f(x)|\le C\,r^{2s-\f{n}{q}}\biggl(\int_{B_r(x)}|f(y)|^q\,dy\biggr)^{\f{1}{q}}.
\end{equation}
Then it follows from (6.4) and changing the order of integrations that
\begin{equation}\begin{split}
\|V(\cS^1_V f)\|^q_{L^q(\BR^n)}&\le C\int_{\BR^n}\biggl(\int_{B_{\f{1}{\fm_V(x)}}(x)}|f(y)|^q\,dy\biggr)\f{V^q(x)}{\fm_V^{2sq-n}(x)}\,dx\\
&=C\int_{\BR^n}|f(y)|^q\biggl(\int_{B_{\f{1}{\fm_V(y)}}(y)}\f{V^q(x)}{\fm_V^{2sq-n}(x)}\,dx\biggr)\,dy.
\end{split}\end{equation}
By Lemma 3.1 and Lemma 3.3, we obtain that
\begin{equation}\begin{split}
\int_{B_{\f{1}{\fm_V(y)}}(y)}\f{V^q(x)}{\fm_V^{2sq-n}(x)}\,dx
&\le\f{C}{\fm_V^{2sq-n}(y)}\int_{B_{\f{1}{\fm_V(y)}}(y)}V^q(x)\,dx\\
&\le\f{C}{\fm_V^{2sq}(y)}\biggl(\fm_V^n(y)\int_{B_{\f{1}{\fm_V(y)}}(y)}V(x)\,dx\biggr)^q\\&\le C.
\end{split}\end{equation}
From (3.2), (6.5) and (6.6), we have that
\begin{equation}\|V(\cS^1_V f)\|_{L^q(\BR^n)}\le C\,\|f\|_{L^q(\BR^n)}\,\,\text{ for $q>\f{n}{2s}$.}
\end{equation}
Using H\"older's inequality as in (6.4) and applying (3.2), Lemma 3.1, Lemma 3.3 and changing the order of integrations, we obtain that
\begin{equation}\begin{split}
\|V(\cS^1_V f)\|_{L^1(\BR^n)}&\le C\int_{\BR^n}|f(y)|\biggl(\int_{B_{\f{1}{\fm_V(y)}}(y)}\f{V(x)}{|x-y|^{n-2s}}\,dx\biggr)\,dy\\
&\le C\int_{\BR^n}|f(y)|\,\f{1}{\fm^{2s}_V(y)}\biggl(\fm^{n}_V(y)\int_{B_{\f{1}{\fm_V(y)}}(y)}V^q(x)\,dx\biggr)^{\f{1}{q}}\,dy\\
&\le C\int_{\BR^n}|f(y)|\biggl(\fm^{n-2s}_V(y)\int_{B_{\f{1}{\fm_V(y)}}(y)}V(x)\,dx\biggr)\,dy\\
&\le C\,\|f\|_{L^1(\BR^n)}.
\end{split}\end{equation}
From standard interpolation argument between the estimates (6.7) and (6.8), we have that
\begin{equation}\|(M_V\circ\cS^1_V) f\|_{L^p(\BR^n)}\le C\,\|f\|_{L^p(\BR^n)}\,\,\text{ for any $p$ with $1\le p\le q$.}
\end{equation}

To deal with $\cS^2_V f(x)$, we note that Corollary 5.8 and H\"older's inequality yield 
\begin{equation*}\begin{split}
|\cS^2_V f(x)|&\le C_N\int_{B^c_r(x)}\f{|f(y)|}{(1+|x-y|\fm_V(x))^N |x-y|^{n-2s}}\,dy\\
&\le C_N\,r^{2s(1-1/p)}\biggl(\int_{B^c_r(x)}\f{|f(y)|^p}{(1+|x-y|\fm_V(x))^N |x-y|^{n-2s}}\,dy\biggr)^{\f{1}{p}}
\end{split}\end{equation*} for $1\le p\le q$ and $r=1/\fm_V(x)$, provided that $N>2s$. Thus we have that
\begin{equation}\begin{split}
&\|V(\cS^2_V f)\|^p_{L^p(\BR^n)}\\
&\quad\le C_N\int_{\BR^n}|f(y)|^p\biggl(\int_{B^c_{\f{1}{\fm_V(y)}}(y)}\f{\fm^{2s(1-p)}_V(x)\,V^p(x)\,dx}{(1+|x-y|\fm_V(x))^N|x-y|^{n-2s}}\biggr)\,dy.
\end{split}\end{equation}
If we set $N_1=\f{N-2s(p-1)d_0}{d_0+1}$ and $S_k=B_{\f{2^k}{\fm_V(y)}}(y)\s B_{\f{2^{k-1}}{\fm_V(y)}}(y)$ for each $k\in\BN$, then it follows from (c) of Lemma 3.1, (3.2), (3.10) and Lemma 3.3 that
\begin{equation}\begin{split}
&\int_{B^c_{\f{1}{\fm_V(y)}}(y)}\f{\fm^{2s(1-p)}_V(x)\,V^p(x)\,dx}{(1+|x-y|\fm_V(x))^N|x-y|^{n-2s}}\\
&\qquad\le\int_{B^c_{\f{1}{\fm_V(y)}}(y)}\f{V^p(x)\,dx}{\fm^{2s(p-1)}_V(y)(1+|x-y|\fm_V(y))^{N_1}|x-y|^{n-2s}}\\
&\qquad\le\f{C}{\fm^{2sp}_V(y)}\sum_{k=1}^{\iy}\f{2^{2s(k-1)}}{(1+2^{k-1})^{N_1}}\,\f{1}{\bigl|B_{\f{2^k}{\fm_V(y)}}(y)\bigr|}\int_{S_k}V^p(x)\,dx\\
&\qquad\le\f{C}{\fm^{2sp}_V(y)}\sum_{k=1}^{\iy}\f{2^{2s(k-1)}}{(1+2^{k-1})^{N_1}}\biggl(\f{1}{\bigl|B_{\f{2^k}{\fm_V(y)}}(y)\bigr|}\int_{B_{\f{2^k}{\fm_V(y)}}(y)}V(x)\,dx\biggr)^p\\
&\qquad\le\f{C}{\fm^{2sp}_V(y)}\sum_{k=1}^{\iy}\f{2^{2s(k-1)}c_1^k}{(1+2^{k-1})^{N_1}}\biggl(\f{1}{\bigl|B_{\f{1}{\fm_V(y)}}(y)\bigr|}\int_{B_{\f{1}{\fm_V(y)}}(y)}V(x)\,dx\biggr)^p\\
&\qquad\le C\sum_{k=1}^{\iy}2^{-(N_1-2s-\log_2 c_1)}\biggl(\fm^{n-2s}_V(y)\int_{B_{\f{1}{\fm_V(y)}}(y)}V(x)\,dx\biggr)^p\\
&\qquad\le C,\,\,\,\,\text{ provided that $N>2(q-1)$ is chosen sufficiently large.}
\end{split}\end{equation}
From (6.10) and (6.11), we thus conclude that
\begin{equation}\|(M_V\circ\cS^2_V) f\|^p_{L^p(\BR^n)}\le C\,\|f\|_{L^p(\BR^n)}\,\,\,\,\text{ for any $p$ with $1\le p\le q$.}
\end{equation}
From the definition of $\cS_V$, we observe that $$(L_K\circ\cS_V)f=f-(M_V\circ\cS_V)f$$ for $f\in L^p(\BR^n)$ with $1\le p\le q$. Therefore the required result can easily be obtained from (6.9) and (6.12). \qed

\,\,\,{\bf [Proof of Theorem 1.4.]}
(a) For $\theta\in[0,n)$, let $\fM_{\theta}g$ be the fractional maximal operator defined by
$$\fM_{\theta} g(x)=\sup_{B\ni x}\f{1}{|B|^{1-\f{\theta}{n}}}\int_B|g(y)|\,dy,$$
where the supremum is taken over every ball $B$ containing $x$. Then it is well-known in standard harmonic analysis \cite{St} that there is a constant $C=C(n,p,q)>0$ such that
\begin{equation}\|\fM_{\theta}g\|_{L^q(\BR^n)}\le C\,\|g\|_{L^p(\BR^n)}
\end{equation} for any $p,q$ with $1<p\le q<\iy$ and $\theta=\ds n\biggl(\f{1}{p}-\f{1}{q}\biggr)$. Thus the proof of (a) of Theorem 1.4 can easily derived from the following lemma.

\begin{lemma}  Let $s\in(0,1)$, $n\ge 2$ and $\theta\in[0,2s)$, and let $\,V\in\rh^\tau$ be nonnegative for $\tau>\f{n}{2s}.$ Then there is a constant $C=C(n,s,\lambda,\theta)>0$ such that
$$\bigl|(M_W\circ\cS_V) f(x)\bigr|\le C\,\fM_{\theta}f(x)$$ for any $x\in\BR^n$.
\end{lemma} 

\pf Take any $x\in\BR^n$ and set $\fm_V(x)=1/\rho$ for $V\in\rh^\tau$ with $\tau>\f{n}{2s}$ and $s\in(0,1)$. If $0\le\theta<2s$, then it follows from (6.2) and Corollary 5.8 that
\begin{equation*}\begin{split}
&|(M_W\circ\cS_V)f(x)|=\bigl|\fm_V^{2s-\theta}(x)\,\cS_V f(x)\bigr|\\
&\qquad\qquad\le C\,\fm_V^{2s-\theta}(x)\int_{\BR^n}\f{|f(y)|}{\bigl(1+|x-y|\,\fm_V(x)\bigr)^{2(2s-\theta)}\,|x-y|^{n-2s}}\,dy\\
&\qquad\qquad\le C\sum_{k=-\iy}^{\iy}\int_{A^k_{\rho}(x)}\f{|f(y)|}{\rho^{2s-\theta}(1+2^{k-1})^{2(2s-\theta)}\,(2^{k-1}\rho)^{n-2s}}\,dy\\
&\qquad\qquad\le C\sum_{k=-\iy}^{\iy}\f{(2^k)^{n-\theta}}{(1+2^{k-1})^{2(2s-\theta)}\,(2^{k-1})^{n-2s}}\biggl(\f{1}{(2^k\rho)^{n-\theta}}\int_{B_{2^k\rho}(x)}|f(y)|\,dy\biggr)\\
&\qquad\qquad\le C\,\fM_{\theta}f(x)\sum_{k=-\iy}^{\iy}\f{(2^k)^{n-\theta}}{(1+2^{k-1})^{2(2s-\theta)}\,(2^{k-1})^{n-2s}}\\
&\qquad\qquad\le C\,\fM_{\theta}f(x)\biggl(\,\,\sum_{k=1}^{\iy}\f{1}{2^{k(2s-\theta)}}+\sum_{k=0}^\iy 2^{-k(2s-\theta)}\biggr)\le C\,\fM_{\theta}f(x),
\end{split}\end{equation*}
where $A^k_{\rho}(x)=B_{2^k\rho}(x)\s B_{2^{k-1}\rho}(x)$. Hence we are done. \qed

\,\,\,(b) $\&$ (c) We note that
\begin{equation}\begin{split}\f{1}{|B|^{1-\f{\theta}{n}}}\int_B|g(y)|\,dy&=\f{|B_1(0)|^{\f{\theta}{n}-1}}{(|B|/|B_1(0)|)^{1-\f{\theta}{n}}}\int_{|y-x|<(\f{|B|}{|B_1(0)|})^{1/n}}|g(y)|\,dy\\
&\le|B_1(0)|^{\f{\theta}{n}-1}\int_{\BR^n}\f{|g(y)|}{|x-y|^{n-\theta}}\,dy
\end{split}\end{equation} for any ball $B\subset\BR^n$ with center $x\in\BR^n$. This imples that the fractional maximal operator $\fM_{\theta}$ is dominated by the Riesz potential $\cI_{\theta}$, i.e.
\begin{equation}\fM_{\theta}g(x)\le C_{n,\theta}\,\,\cI_{\theta}|g|(x)=C_{n,\theta}\,\,|g|*I_\theta(x)
\end{equation} where $C_{n,\theta}>0$ is certain constant depending only on $n$ and $\theta\in[0,n)$ and $I_\theta$ is the Riesz kernel given by $I_\theta(y)=|y|^{-n+\theta}$. Then it is easy to check that 
\begin{equation}\|I_\theta\|_{L^{\f{n}{n-\theta},\iy}(\BR^n)}<\iy\,\,\text{ for any $\theta\in[0,n)$.}
\end{equation}
Also, it is well-known \cite{St} that there is a universal constant $C(n,s)>0$ such that
\begin{equation}\|I_\theta g\|_{L^{q,\iy}(\BR^n)}\le C(n,s)\|g\|_{L^1(\BR^n)}
\end{equation}
for any $q\in(1,\f{n}{n-2s})$. Thus we can easily derive (b) from (6.15) and (6.17).
Finally, we can easily obtain (c) from (6.16) and Proposition 7.1 below, because
$$\f{1}{p}-\f{1}{q}=1-\f{1}{r}=\f{2s}{n}$$
when $r=\f{n}{n-2s}$. Hence we complete the proof. \qed

\section{Appendix}

In order to obtain the mapping properties of $M_W\circ\cS_V$ on the boundary of the trapezoidal area in Figure 1, we need the following estimates whose proof is self-contained. 

\begin{prop} If $p\in[1,\iy)$ and $q,r\in(1,\iy)$ satisfy that
\begin{equation}\f{1}{p}-\f{1}{q}=1-\f{1}{r},
\end{equation} then there is a constant $C=C(p,q,r)>0$ such that
\begin{equation*}\|g*h\|_{L^{q,\iy}(\BR^n)}\le C\,\|h\|_{L^{r,\iy}(\BR^n)}\|g\|_{L^p(\BR^n)}
\end{equation*} for any $g\in L^p(\BR^n)$ and $h\in L^{r,\iy}(\BR^n)$. Moreover, $C=\cO((r-1)^{-p/q})$ as $r\to 1^-$.
\end{prop}

\pf For $N>0$ and $h\in L^{r,\iy}(\BR^n)$, we denote by 
\begin{equation*}L_N=\{y\in\BR^n:|h(y)|\le N\}\,\,\text{ and }\,\,U_N=\{y\in\BR^n:|h(y)|> N\}.
\end{equation*}
If we set $h_1=h\mathbbm{1}_{L_N}$ and $h_2=h\mathbbm{1}_{U_N}$, then we have that
\begin{equation}\begin{split}\om_{h_1}(\gm)&=[\om_h(\gm)-\om_h(N)]\mathbbm{1}_{(0,N)}(\gm),\\
\om_{h_2}(\gm)&=\om_h(N)\mathbbm{1}_{(0,N]}(\gm)+\om_h(\gm)\mathbbm{1}_{(N,\iy)}(\gm)\,\,\text{ for any $\gm>0$.}
\end{split}\end{equation}
It is easy to check that
\begin{equation}\om_{g*h}(\gm)\le\om_{g*h_1}(\gm/2)+\om_{g*h_2}(\gm/2)\,\,\text{ for any $\gm>0$.}
\end{equation}
From (7.1), we see that $1<r<p'$ where $p'$ is the dual exponent of $p$. If $p'<\iy$, then by simple calculation and (7.2) we have that
\begin{equation}\begin{split}
\int_{\BR^n}|h_1(y)|^{p'}\,dy&=p'\int_0^{\iy}\gm^{p'-1}\om_{h_1}(\gm)\,d\gm\\
&=p'\int_0^N\gm^{p'-1}[\om_h(\gm)-\om_h(N)]\,d\gm\\
&\le p'\int_0^N\gm^{p'-1-r}\|h\|^r_{L^{r,\iy}(\BR^n)}\,d\gm
-p'\int_0^N\gm^{p'-1}\om_h(N)\,d\gm\\
&=\f{p' N^{p'-r}}{p'-r}\|h\|^r_{L^{r,\iy}(\BR^n)}-N^{p'}\om_h(N)\le\f{p' N^{p'-r}}{p'-r}\|h\|^r_{L^{r,\iy}(\BR^n)},
\end{split}\end{equation}
and so it follows from H\"older's inequality and (7.4) that
\begin{equation}\begin{split}\|g*h_1\|_{L^{\iy}(\BR^n)}&\le\|g\|_{L^p(\BR^n)}\|h_1\|_{L^{p'}(\BR^n)}\\&\le\|g\|_{L^p(\BR^n)}\biggl(\f{p' N^{p'-r}}{p'-r}\|h\|^r_{L^{r,\iy}(\BR^n)}\biggr)^{1/p'}.
\end{split}\end{equation}
If $p'=\iy$, then by applying H\"older's inequality again we obtain that
\begin{equation}\|g*h_1\|_{L^{\iy}(\BR^n)}\le\|g\|_{L^p(\BR^n)}\,N.
\end{equation}
Fix any $\gm>0$ and choose some $N$ so that
\begin{equation}\f{\gm}{2\|g\|_{L^p(\BR^n)}}=\begin{cases}\biggl(\ds\f{p' N^{p'-r}}{p'-r}\|h\|^r_{L^{r,\iy}(\BR^n)}\biggr)^{1/p'} &\text{ if $p'<\iy$, }\\
\qquad\qquad\qquad N &\text{ if $p'=\iy$.}\end{cases}
\end{equation}

On the other hand, by simple calculation and (7.2), we have that
\begin{equation}\begin{split}
\int_{\BR^n}|h_2(y)|\,dy&=\int_0^{\iy}\om_{h_2}(\gm)\,d\gm=\int_0^N\om_h(N)\,d\gm+\int_N^{\iy}\om_h(\gm)\,d\gm\\
&\le N\,\om_h(N)+\int_N^{\iy}\gm^{-r}\,\|h\|^r_{L^{r,\iy}(\BR^n)}\,d\gm\\
&\le N^{1-r}\|h\|^r_{L^{r,\iy}(\BR^n)}+\f{N^{1-r}}{r-1}\|h\|^r_{L^{r,\iy}(\BR^n)}\\
&=\f{r N^{1-r}}{r-1}\|h\|^r_{L^{r,\iy}(\BR^n)}, 
\end{split}\end{equation}
and thus it follows from Young's inequality and (7.8) that
\begin{equation}\|g*h_2\|_{L^p(\BR^n)}\le\|g\|_{L^p(\BR^n)}\|h_2\|_{L^1(\BR^n)}\le\f{r N^{1-r}}{r-1}\|h\|^r_{L^{r,\iy}(\BR^n)}\|g\|_{L^p(\BR^n)}.
\end{equation}
By (7.6), (7.7), (7.9) and Chebychev's inequality, we conclude that
\begin{equation}\begin{split}
\om_{g*h}(\gm)&\le\om_{g*h_2}(\gm/2)\le\f{2^p}{\gm^p}\,\|g*h_2\|^p_{L^p(\BR^n)}\\
&\le\f{2^p}{\gm^p}\,\biggl(\f{r N^{1-r}}{r-1}\|h\|^r_{L^{r,\iy}(\BR^n)}\|g\|_{L^p(\BR^n)}\biggr)^p\\
&=\f{C^q}{\gm^q}\|h\|^q_{L^{r,\iy}(\BR^n)}\|g\|^q_{L^p(\BR^n)}.
\end{split}\end{equation}
Therefore the required result can be obtained from (7.10). \qed


\end{document}